\documentclass[reqno]{amsart}
\usepackage{amssymb,bbm,graphicx,mathtools,verbatim}
\usepackage[colorlinks]{hyperref}
\usepackage[text={360pt,594pt},centering]{geometry}
\newtheorem{thm}{Theorem}[section]
\newtheorem{cor}[thm]{Corollary}
\newtheorem{lem}[thm]{Lemma}

\theoremstyle{definition}

\newtheorem{exm}[thm]{Example}
\numberwithin{equation}{section}

\DeclarePairedDelimiter{\norm}{\lVert}{\rVert}

\newcommand{\ov}{\overline}
\renewcommand{\Re}{\operatorname{Re}}
\renewcommand{\Im}{\operatorname{Im}}

\newcommand{\e}{\mathrm{e}}
\newcommand{\iu}{\mathrm{i}}

\newcommand{\sm}[1]{\big(\begin{smallmatrix}#1\end{smallmatrix}\big)}

\newcommand{\bb}[1]{{\mathbb{#1}}}
\newcommand{\mc}[1]{{\mathcal{#1}}}
\newcommand{\id}{\mathbbm 1}

\newcommand{\ol}{{\overline{\lambda}}}
\newcommand{\la}{\lambda}

\newcommand{\sas}{\;\;\text{and}\;\;}

\begin{document}
\title[The limit-point/limit-circle classification]{The limit-point/limit-circle classification for ordinary differential equations with distributional coefficients}
\author{Varun Bhardwaj and Rudi Weikard}
\address{Department of Mathematics, University of Alabama at Birmingham, Birmingham, AL 35226-1170, USA}
\email{varun01@uab.edu, weikard@uab.edu}
\date{\today}


\begin{abstract}
We investigate the limit-point/limit-circle classification for the differential equation $Ju'+qu=\la wu$ where $J=\sm{0&-1\\ 1&0}$ and $q$ and $w$ are matrices whose entries are distributions of order zero with $q$ Hermitian and $w$ non-negative.
We identify the situations when the classical alternative of Weyl works and when it fails.
\end{abstract}
\maketitle

\section{Introduction}
In his ingenious work \cite{41.0343.01} of 1910 Weyl showed how to determine whether a Sturm-Liouville equation with one singular and one regular endpoint had one or two linearly independent square integrable solutions by approximating it with regular problems.
An answer to that question is important in the context of the extension theory for the associated symmetric operator.
Specifically, posing the equation $-(py')'+vy=\la ry$ on the interval $(0,\infty)$ with $0$ a regular endpoint and $\la$ a parameter in $\bb C\setminus\bb R$, Weyl introduced an auxiliary boundary condition at a finite point $c$, namely $\cos(\beta)y(c)+\sin(\beta)(py')(c)=0$.
Among the solutions of this boundary value problem there is one of the form $\varphi+m(\beta)\psi$ where $\varphi$ and $\psi$ are solutions satisfying appropriate initial condition at $0$.
As $\beta$ varies in $[0,\pi)$ the coefficient $m(\beta)$ describes a circle in the complex plane.
Now, as $c$ approaches infinity the disks bounded by these circles are nested and shrink either to a point (the limit-point case) or to a disk (the limit-circle case).
In the former case only one linearly independent solution of the original problem is in the weighted Hilbert space $L^2(r)$.
In the latter case all solutions are in $L^2(r)$.
Accordingly one can pose either one boundary condition (at $0$) or two boundary conditions (separated or coupled) in order to arrive at self-adjoint realizations of the differential expression $(-(py')'+qy)/r$.
This relationship is known as Weyl's alternative.

In a series of papers from around 1940 (later collected in \cite{MR0176151}) Titchmarsh investigated this problem using extensively the tools of Complex Analysis.
He showed that the coefficient $m$, now called the Weyl-Titchmarsh coefficient, must satisfy $\Im(m)/\Im(\la)>0$.
Since, as a function of $\la$, it is also analytic away from the real axis, it is a Nevanlinna function and as such has an integral representation with respect to a positive measure called the spectral measure, an indispensable tool for the study of the equation.
For further information on the history of Sturm-Liouville theory we refer to L\"utzen~\cite{MR745152}.
Everitt \cite{MR2145077} has a shorter account focusing on the contributions of Weyl and Titchmarsh.

The goal of this paper is to show that Weyl's ideas largely work for the case of the differential equation
$$Ju'+qu=\la w u,$$
posed on the interval $(0,b)$ (with $0<b\leq\infty$),
where $J=\sm{0&-1\\ 1&0}$ while $q$ and $w$ are $2\times 2$-matrices whose entries are distributions of order $0$ on $(0,b)$.
While we shall not require that $q$ and $w$ are real we do require that $q$ is hermitian and that $w$ is non-negative (but not identically $0$) allowing for complex distributions in the off-diagonal entries of $q$ and $w$.
We still require that $0$ is a regular point since one can deal with two singular endpoints in exactly the same way as in the classical case.

The Sturm-Liouville equation is subsumed under our approach.
In fact the more general equation $-(p(y'+sy))'+sp(y'+sy)+vy=\la r y$, treated by Eckhardt et al. in \cite{MR3046408} corresponds to the system $Ju'+qu=\la wu$ upon setting $u=\sm{y\\ p(y'+sy)}$, $q=\sm{v&s\\ s&-1/p}$, and $w=\sm{r&0\\ 0&0}$ where $r$, $q$, $s$, and $1/p$ are real, locally integrable, and integrable near $0$.

There is an excellent account of Weyl's circle of ideas in Coddington and Levinson~\cite{MR0069338}.
Our approach to the problem emulates theirs.
However, before we get started we will present some background information.

Antiderivatives of distributions of order $0$ on $(0,b)$ may be represented by functions of locally bounded variation.
Thus, if $s$ is a distribution of order $0$ with antiderivative $S$ we have, for all test functions $\phi$, that
$$s(\phi)=\int_{(0,b)} \phi\;dS$$
where $dS$ is the local\footnote{In general $dS$ is a measure only on compact subsets of $(0,b)$, hence the modifier ``local''.} Lebesgue-Stieltjes measure generated by $S$.
Thus, in the classical case, the local (matrix-valued) measures $q$ and $w$ are locally absolutely continuous with respect to Lebesgue measure.
By saying that $0$ is a regular point we mean that $q$ and $w$ are (finite) measures near $0$.

Solutions $u$ of the equation are to be thought of as functions of locally bounded variation so that $u'$ is also a distribution of order $0$.
We refer to \cite{MR4047968} for some background about ordinary differential equations with distributional coefficients and further information on the relationship between distributions of order $0$ and local measures.

$w$, due to it being non-negative, is, in fact, a matrix-valued measure.
Therefore $\int_{(0,b)} f^*wg$ is a semi-scalar product on the space $\mc L^2(w)$ of functions $f$ satisfying $\int_{(0,b)} f^*wf<\infty$.
The space $L^2(w)$, defined as the quotient of $\mc L^2(w)$ and the space of the functions of norm $0$, is a Hilbert space.

When the local measures $q$ or $w$ have atomic parts there is an issue with existence and uniqueness of solutions of initial value problems.
We define the matrices\footnote{In the following we use the notation $u^+(x)=\lim_{t\downarrow x}u(t)$ and $u^-(x)=\lim_{t\uparrow x}u(t)$.} $\Delta_q(x)=Q^+(x)-Q^-(x)=q(\{x\})$ and $\Delta_w(x)=W^+(x)-W^-(x)=w(\{x\})$ where $Q$ and $W$ are, respectively, antiderivatives of $q$ and $w$ and therefore functions of locally bounded variation.
We also set
$$B_{\pm}(x,\lambda)=J\pm \frac{1}{2}(\Delta_q(x)-\lambda \Delta_w(x)).$$
According to Theorem 2.2 of \cite{MR4047968} a unique balanced\footnote{We use the notation $u^\#(x)=(u^-(x)+u^+(x))/2$ and call a function $u$ of locally bounded variation balanced, if $u=u^\#$.} solution of $Ju'+qu=\la w u$ satisfying the initial condition $u(x_0)=u_0\in\bb C^2$ exists in the interval $(c,d)\subset(0,b)$, if $x_0\in(c,d)$ and $B_\pm(x,\la)$ are invertible for all $x\in(c,d)$.
If $c=0$, the initial condition can also be posed at $0$ due to it being a regular endpoint.
Since $B_+(x,\la)u^+(x)=B_-(x,\la)u^-(x)$ is equivalent to $Ju'(\{x\})+(q(\{x\})-\la w(\{x\}))u(x)=0$,
it is clear that the lack of invertibility of either $B_+(x,\la)$ or $B_-(x,\la)$ causes problems.
We introduce the set $\Lambda=\{\la\in\bb C:\exists x\in(0,b): \det(B_+(x,\la))\det(B_-(x,\la))=0\}$ of those $\la$ where invertibility lacks for one of the matrices $B_\pm(x,\la)$ for at least one $x\in(0,b)$ and require throughout the paper that $\Lambda\neq\bb C$.
We call a point $x$ a bad point for $\la$ if one of $B_\pm(x,\la)$ is not invertible.
Note that $B_\pm(x,\la)=-B_\mp(x,\ol)^*$.
Therefore the set $\Lambda$ is symmetric with respect to the real axis.
Since $Q$ and $W$ are of locally bounded variation, they have only countably many jumps.
However, these jumps may form a dense subset of $(0,b)$.
A point $x$ at which each entry of $Q$ and $W$ is continuous will be called a point of continuity below.

It may happen\footnote{This is so even for constant coefficients. It is not a consequence of distributional coefficients.} that the space of solutions of the differential equation $Ju'+qu=0$ which have norm $0$, a space we denote by $\mc L_0$, is not trivial.
If $u\in\mc L_0$, it also satisfies $Ju'+qu=\la wu$ for any $\la\in\bb C$ since $\norm u=0$ if and only if $wu=0$ as a distribution (or measure).
We say that the definiteness condition holds if $\mathcal{L}_0=\{0\}$.
Since we exclude the trivial case $w=0$ and since $\Lambda\neq\bb C$ we know that $\mathcal{L}_0$ is either trivial or one-dimensional.

It follows from the theory of symmetric relations that the number of linearly independent solutions of $Ju'+qu=\la wu$ of positive finite norm does not change as $\la$ varies in either the upper or the lower half of the complex plane.
These numbers are called deficiency indices and are denoted by $n_+$ (upper half-plane) and $n_-$ (lower half-plane).
Since a solution of norm $0$ for some $\la$ is also a solution for any other $\la$, the number of linearly independent solutions of $Ju'+qu=\la wu$ of finite norm also depends only on whether $\la$ is in the upper or the lower half plane.

To finish the introduction we will summarize the most important results obtained below.
There is a function $\tau$ such that
$$\tau(x,\la)=\frac{u^-(x)}{\ov{v^-(x)}}$$
whenever $Ju'+qu=\la wu$, $Jv'+qv=\ol wv$ and $\ov{v(0)}=u(0)$.
This function satisfies $\tau(x,\la)\ov{\tau(x,\ol)}=1$.
If $q$ and $w$ are real, then $\tau(x,\la)=1$.
The deficiency indices $n_+$ and $n_-$ are identical if neither $\tau(x,\la)$ nor $\tau(x,\ol)$ tends to $0$ as $x$ tends to $b$.
We emphasize that this condition is sufficient but not necessary as we show by an example.

If $\la$ is not in $\Lambda$ and if the definiteness condition holds we have (eventually) nesting Weyl disks.
The disks associated with $\ol$ are reflections across the real line of those associated with $\la$.
If $n_+=n_-$ we have, in resemblance of the classical Weyl alternative, $n_\pm=1$ or $n_\pm=2$ precisely when the disks shrink to a point or to a disk, respectively.
If the definiteness condition does not hold, we may have nested half-planes instead of disks.
Again, if $n_+=n_-$, we have that $n_\pm=0$ or $n_\pm=1$ precisely when the half planes or disks shrink to a point, possibly $\infty$, or to a half-plane or disk, respectively.

Finally we show that it is important to choose $\la$ outside $\Lambda$, since we always have the limit-point case at $b$ otherwise, even if $b$ is a regular endpoint.

\section{A sufficient condition for equal deficiency indices}
\begin{lem}
Suppose $\la\not\in\Lambda$ and $n\mapsto x_n$ is an enumeration of all points $x$ in $(0,b)$ where $\Delta_q(x)$ or $\Delta_w(x)$ is different from $0$.
Then
\[P(x,\lambda)=\prod_{x_n \in (0,x)}  \frac{\det B_-(x_n,\lambda)}{\det B_+(x_n,\lambda)}\]
exists for every $x\in [0,b)$, is independent of the order of the factors, and is different from $0$.
\end{lem}

\begin{proof}
Fix $x\in[0,b)$.
The cases when $x=0$ or when the set $\{x_n:x_n \in (0,x)\}$ is empty, are trivial if we define the empty product to be $1$.
Note that
\begin{equation}\label{detbpm}
\det B_-(x_n,\la)-\det B_+(x_n,\la)=2\iu(\Im\Delta_q(x_n)_{12}-\la \Im\Delta_w(x_n)_{12}).
\end{equation}
Since $q$ and $w$ are finite measures on $(0,x)$ we have that the entries of $\Delta_q(x_n)$ and $\Delta_w(x_n)$ are absolutely summable, if we consider only $x_n\in(0,x)$.
It follows that $\det B_-(x_n,\la)-\det B_+(x_n,\la)$ is absolutely summable and that $\det B_+(x_n,\la)$ is bounded away from $0$.
Therefore $\sum_{x_n \in (0,x)}\left|\frac{\det B_-(x_n,\lambda)}{\det B_+(x_n,\lambda)}-1\right|$ is convergent which, as is well-known, implies our claim.
\end{proof}

For $\la\not\in\Lambda$ we define the function
\[\tau(x,\lambda)=P(x,\lambda) \exp\left(2\iu \int_{(0,x)} (\Im(q_{c,12})-\lambda \Im(w_{c,12})) \right),\]
where $q_{c,12}$ and $w_{c,12}$ denote the continuous parts of the (local) measures $q_{12}$ and $w_{12}$, respectively.
Note that $\tau(\cdot,\la)$ is left-continuous and that $\tau^+(x,\la)=\tau(x,\la)\det B_-(x,\la)/\det B_+(x,\la)$.
Also, $\tau(x,\la)=1$ when $q$ and $w$ are real.
Since $B_\pm(x,\la)=-B_\mp(x,\ol)^*$ we have $\det B_\pm(x,\la)=\ov{\det B_\mp(x,\ol)}$.
This implies immediately that
\begin{equation}\label{tau}
\tau(x,\la)\ov{\tau(x,\ol)}=1.
\end{equation}

\begin{lem}\label{L_Temp_01}
Suppose $\la\not\in\Lambda$.
Let $u$ be a balanced solution of $Ju'+qu=\la wu$ and define $v^-(x)=\tau(x,\ol) \overline{u^-(x)}$.
Then $v=v^\#$, the balanced version of $v^-$, satisfies $Jv'+qv=\ol wv$ and $v(0)=\overline{u(0)}$.
\end{lem}

\begin{proof}
To prove that $v$ satisfies the differential equation, we will show that the measure $Jv'+qv-\ol wv$ is the zero measure.
We will do this in two steps.
First we will show that it is a continuous measure, i.e. $(Jv'+qv-\ol wv)(\{x\})=0$ for all $x$, and then complete the proof by showing that
$(Jv'+qv-\ol wv)((0,x))=0$.
We know that
\begin{equation}\label{equ}
B_+(x,\la) u^+(x)=B_-(x,\la) u^-(x)
\end{equation}
and we have to show that
\begin{equation}\label{eqv}
B_+(x,\ol)\tau^+(x,\ol)\overline{u^+(x)}=B_-(x,\ol)\tau(x,\ol)\overline{u^-(x)}.
\end{equation}

As mentioned above we have $B_\pm(x,\ol)=-B_\mp(x,\la)^*=-\ov{B_\mp(x,\la)^\top}$.
We therefore get from equation \eqref{equ} that
\[J^{-1}B_-(x,\ol)^\top \ov{u^+(x)}=J^{-1}B_+(x,\ol)^\top \ov{u^-(x)}.\]
Next note that, for any invertible $2\times 2$-matrix $M$ we have $(\det M) M^{-1}=J^{-1}M^\top J$ which leads to
\[\tau^+(x,\ol) B_-(x,\ol)^{-1}J^{-1}\ov{u^+(x)}=\tau(x,\ol) B_+(x,\ol)^{-1}J^{-1}\ov{u^-(x)}.\]
Finally, for any $2\times 2$-matrix $M$ we have $(J+M)J(J-M)=(J-M)J(J+M)$ which implies that $B_-(x,\ol)B_+(x,\ol)^{-1}J^{-1}=J^{-1}B_+(x,\ol)^{-1}B_-(x,\ol)$.
Applying this to the previous identity we get equation \eqref{eqv} and hence that the measure $Jv'+qv-\ol wv$ is continuous.

Next, define the set $A=\{y\in(0,x): \Delta_q(y)=0=\Delta_w(y)\}$.
Since the complement $A^c$ of $A$ is countable we have $\int_{(0,x)} (Jv'+(q-\ol w)v)=\int_A (Jv'+(q-\ol w)v)$.
The measures $v'$ and $(v^-)'$ are identical and, on $A$, we have $\tau^+=\tau^-=\tau$ as well as $u^+=u^-$.
Therefore
$$\int_A (Jv'+(q-\ol w)v)=\int_A \big(\tau'(\cdot,\ol)J\ov u+\tau(\cdot,\ol)(J\ov u'+(q-\ol w)\ov u)\big).$$
Since $q=\ov q-2\iu \Im q_{12}J$ and $w=\ov w-2\iu \Im w_{12}J$ we obtain
$$J\ov u'+(q-\ol w)\ov u=-2\iu (\Im q_{12}-\ol \Im w_{12})J\ov u.$$
This implies that
$$\int_A (Jv'+(q-\ol w)v)=\int_A \big(\tau'(\cdot,\ol)-2\iu(\Im q_{12}-\ol \Im w_{12})\tau(\cdot,\ol)\big)J\ov u.$$
Our definition of $\tau$ shows now that $\int_A (Jv'+(q-\ol w)v)=0$ and thus that $Jv'+qv=\ol wv$.

The claim about the initial condition of $u$ and $v$ follows since $\tau(0,\la)=1$.
\end{proof}

In the following we will have frequent need for Green's formula or Lagrange's identity, a major tool of spectral theory.
It still works under our present hypotheses (see equation (3.3) in \cite{MR4047968}) and states that
\begin{equation}\label{Lagrange}
(v^*Ju)^-(d)-(v^*Ju)^+(c)=\int_{(c,d)}(v^*wf-g^*wu)
\end{equation}
provided that $v,g,u,f\in\mc L^2(w)$, $Jv'+qv=wg$, and $Ju'+qu=wf$ on $(c,d)\subset(0,b)$.

\begin{lem}\label{L.2.3}
Suppose $\la\not\in\Lambda$, $u$ and $v$ satisfy $Ju'+qu=\la wu$ and $Jv'+qv=\ol wv$, respectively, and $v(0)=\ov{u(0)}$.
If $\tau(\cdot,\la)$ is bounded and bounded away from $0$, then $u\in \mc L^2(w)$ if and only if $v\in \mc L^2(w)$.
\end{lem}

\begin{proof}
There is nothing to prove when $\la\in\bb R$.
Therefore we may assume that $\Im(\la)\neq0$.
For a $\bb C^2$-valued function $y$ of locally bounded variation define $\norm y_x^2=\int_{(0,x)} y^*wy$.
Choosing $(u,f)=(v,g)=(u,\la u)$ in Lagrange's identity \eqref{Lagrange} gives
$$\norm u_x^2 = \frac1{2\iu\Im\la} \big((u^*Ju)(x)-(u^*Ju)(0)\big)$$
whenever $x$ is a point of continuity.
Similarly,
$$\norm v_x^2 = \frac{-1}{2\iu\Im\la} \big((v^*Jv)(x)-(v^*Jv)(0)\big).$$
Since the existence and uniqueness theorem for initial value problems is in force when $\la\not\in\Lambda$ and since $x$ is a point of continuity, Lemma \ref{L_Temp_01} shows that $v(x)=\tau(x,\ol)\ov{u(x)}$.
This gives
$$\norm v_x^2 = |\tau(x,\ol)|^2\norm u_x^2
-\frac1{2\iu\Im\la}\left(|\tau(x,\ol)|^2-1\right) \ov{(u^*Ju)(0)}$$
where we used equation \eqref{tau}.

Equation \eqref{tau} shows also that $|\tau(x,\ol)|=1/|\tau(x,\la)|$ which is bounded by our hypothesis.
Thus, if $x\mapsto \norm u_x^2$ is bounded then so is $x\mapsto\norm v_x^2$.
The other  direction is, of course, proved in the same way.
\end{proof}

When $\tau(\cdot,\la)$ is bounded and bounded away from zero we have therefore that all solutions of $Ju'+qu=\la wu$ are in $\mc L^2(w)$ if and only if the same is true for all solutions of $Ju'+qu=\ol wu$.

\begin{thm}
Suppose $\la_0$ and $\ol_0$ are in $\bb C\setminus\Lambda$ and $\la\in\bb C$.
If all solutions of $Ju'+qu=\la_0 wu$ and all solutions of $Ju'+qu=\ov{\la_0} wu$ are in $\mc L^2(w)$, then any solution of $Ju'+qu=\la wu$ is also in $\mc L^2(w)$.
In particular, $n_+=n_-=2-\dim\mc L_0$.
\end{thm}

\begin{proof}
We define the norm $\norm u_{c,d}$ by setting
$$\norm u_{c,d}=\int_{(c,d)} u^*wu$$
whenever $c$ and $d$ are points of continuity.

For $\mu=\la_0$ and $\mu=\ol_0$ let $U(\cdot,\mu)$ be the fundamental matrix for $Ju'+qu=\mu wu$ satisfying $U(0,\mu)=\id$ and denote the first and second column of $U(\cdot,\mu)$ by $\varphi(\cdot,\mu)$ and $\psi(\cdot,\mu)$, respectively.
These are elements of $\mc L^2(w)$ and therefore we may choose $c$ such that $\norm{\varphi(\cdot,\la_0)}_{c,d}$ and $\norm{\psi(\cdot,\la_0)}_{c,d}$ as well as $\norm{\varphi(\cdot,\ol_0)}_{c,d}$ and $\norm{\psi(\cdot,\ol_0)}_{c,d}$ are smaller than $M=1/(2\sqrt{|\la-\la_0|})$ for all $d\in(c,b)$.

Now suppose $u$ satisfies $Ju'+qu=\la wu$ or, equivalently, $Ju'+(q-\la_0w)u=(\la-\la_0)wu$.
Then, by the variation of constants formula (cf. Lemma 3.3 of \cite{MR4047968}), there is a vector $(\alpha,\beta)^\top\in\bb C^2$ such that
$$u(x)=U(x,\la_0)\big(\sm{\alpha\\ \beta} +(\la-\la_0)J^{-1}\int_{(c,x)} U(\cdot,\ov{\la_0})^*wu\big)$$
whenever $x\geq c$ is a point of continuity.

Set $\gamma(x)=\int_{(c,x)}\varphi(\cdot,\ov{\la_0})^*wu$ and $\delta(x)=\int_{(c,x)}\psi(\cdot,\ov{\la_0})^*wu$.
Then we have, by Cauchy-Schwarz, $|\gamma(x)|\leq M \norm u_{c,d}$ and $|\delta(x)|\leq M \norm u_{c,d}$ as long as $x\in(c,d)$.
From Minkowski's inequality we get
$$\norm u_{c,d}\leq (|\alpha|+|\la-\la_0| M \norm u_{c,d})M +(|\beta|+|\la-\la_0| M \norm u_{c,d})M.$$
With our definition of $M$ this yields
$$\norm u_{c,d}\leq 2(|\alpha|+|\beta|)M.$$
Here the right hand side is independent of $d$ and this implies that $u\in\mc L^2(w)$.
\end{proof}

\section{The limit-point/limit-circle classification}
We consider the first-order system
$$Ju'+qu=\lambda wu$$
when $\la$ is a fixed number in $\bb C\setminus\bb R$.

\subsection{No bad points exist for \texorpdfstring{$\la$}{lambda}}\label{gps}
In this section we assume that $B_\pm(x,\la)$ are invertible for all $x\in(0,b)$, i.e., $\la\not\in\Lambda$.
Under this assumption we can uniquely solve any initial value problem with the initial condition posed at $0$.
Accordingly we define the fundamental matrix $U(\cdot,\la)$ to be that solution of $Ju'+qu=\la wu$ which satisfies the initial condition
$U(0,\la)=\sm{\cos\alpha&-\sin\alpha\\ \sin\alpha&\cos\alpha}$ for some fixed $\alpha\in[0,\pi)$.
We also denote the first and second columns of $U(\cdot,\la)$ by $\varphi(\cdot,\la)$ and $\psi(\cdot,\la)$, respectively.
For a fixed $c\in(0,b)$, a point of continuity, we now ask for which numbers $m$ the function $\chi_m(\cdot,\la)=U(\cdot,\la)\sm{1\\ m}$ satisfies the boundary condition $(\cos\beta,\sin\beta)\chi_m(c,\la)=0$ for $\beta\in[0,\pi)$.

Setting $z=\cot(\beta)$, $A=U_{11}(c,\lambda)$, $B=U_{21}(c,\lambda)$, $C=U_{12}(c,\lambda)$, and $D=U_{22}(c,\lambda)$, this is equivalent to the requirement
$$(z,1)\begin{pmatrix}A&C\\ B&D\end{pmatrix}\begin{pmatrix}1\\ m\end{pmatrix}=0$$
which yields the Möbius transform\footnote{$U(c,\la)=\sm{A&C\\ B&D}$ is invertible for every $c\in(a,b)$.}
\[m=-\frac{Az+B}{Cz+D}\]
whose inverse is
\[z=-\frac{Dm+B}{Cm+A}.\]
The Möbius transform $m(z)$ maps the real line\footnote{We are a little imprecise here (and elsewhere), it is actually the one-point compactification of the real line we have in mind.} either onto a circle or else a straight line.
For $z\in\bb R$ it follows from $z=\ov z$ that $m$ must satisfy the equation
\begin{equation}\label{col}
(C\ov{D}-\ov{C}D)|m|^2+(A\ov{D}-B\ov{C})\ov m+(\ov{B}C-\ov{A}D)m+A\ov{B}-\ov{A}B=0.
\end{equation}

If $C\overline{D}-\overline{C}D\neq0$ equation \eqref{col} implies that $|m-\tilde{m}|^2=r^2$ where
$$\tilde{m}=\frac{B \overline{C}-A\overline{D}}{C\overline{D}-\overline{C}D} \sas r=\left|\frac{AD-BC}{C\overline{D}-\overline{C}D} \right|.$$
Thus $m$ describes a circle centered at $\tilde{m}$ with radius $r$.
Of course, $\tilde m$ and $r$ depend on $c$ and $\la$.
We denote the disk $\{m\in\bb C:|m-\tilde m|\leq r\}$ by $D(c,\la)$.

If $C\overline{D}-\overline{C}D=0$ we learn from \eqref{col} that
$$\varrho \ov m-\ov\varrho m+\sigma=2\iu\Im(\varrho)\Re(m) -2\iu\Re(\varrho)\Im(m)+\sigma=0$$
where $\varrho=A\ov D-B\ov C$ and $\sigma=A\ov B-\ov A B\in\iu\bb R$.
Thus $m$ describes a straight line with slope $\Im(\varrho)/\Re(\varrho)$.
We will later see that $\varrho=1$, so that the line is horizontal.

Next note that the expression on the left of \eqref{col} equals $\chi_m(c,\la)^*J\chi_m(c,\la)$ while
\begin{equation}\label{cdcd}
C\ov{D}-\ov{C}D=\psi(c,\la)^*J\psi(c,\la).
\end{equation}
Let us introduce the norm $\norm{u}_c$ by setting $\norm{u}_c^2=\int_{(0,c)}u^*wu$ and use Lagrange's identity \eqref{Lagrange} to bring the norms of $\chi_m(\cdot,\la)$ and $\psi(\cdot,\la)$ into the game.
Since $\chi_m(0,\lambda)^*J\chi_m(0,\lambda)=-2\iu\Im(m)$ Lagrange's identity gives
\[\chi_m(c,\lambda)^*J\chi_m(c,\lambda)+2\iu\Im(m)=2\iu \Im(\lambda)\norm{\chi_m(\cdot,\lambda)}_c^2\]
upon choosing $(u,f)=(v,g)=(\chi_m(\cdot,\lambda),\lambda \chi_m(\cdot,\lambda))$.
Similarly,
\begin{equation}\label{lagforpsi}
\psi(c,\lambda)^*J\psi(c,\lambda)=2\iu \Im(\lambda)\norm{\psi(\cdot,\lambda)}_c^2
\end{equation}
using $(u,f)=(v,g)=(\psi(\cdot,\lambda),\lambda \psi(\cdot,\lambda))$ and $\psi(0,\lambda)^*J\psi(0,\lambda)=0$.

When $\norm{\psi(\cdot,\lambda)}_c\neq0$, i.e., when $m(\bb R)$ is a circle, equation \eqref{col} becomes
$$\norm{\chi_m(\cdot,\lambda)}_c^2-\frac{\Im(m)}{\Im(\lambda)}=\big(|m-\tilde{m}|^2-r^2\big)\norm{\psi(\cdot,\lambda)}_c^2=0.$$
Note that the left equation here holds for all $m\in\bb C$.
Therefore, $m\in D(c,\la)$ if and only if
\begin{equation}\label{nchi1}
\norm{\chi_m(\cdot,\lambda)}_c^2\leq\frac{\Im(m)}{\Im(\lambda)}.
\end{equation}

When $\norm{\psi(\cdot,\lambda)}_c=0$ and $m(\bb R)$ is a line, we get instead
$$\Im(\lambda)\norm{\chi_m(\cdot,\lambda)}_c^2-\Im(m)=\Im(\varrho)\Re(m)-\Re(\varrho)\Im(m)+\frac\sigma{2\iu}=0.$$
Since $\norm{\chi_m(\cdot,\lambda)}_c=\norm{\varphi(\cdot,\lambda)}_c$ does not depend on $m$ and since the first of these equations holds for all $m\in\bb C$, we see that $\varrho=1$, i.e., the line $m(\bb R)$ is horizontal and is given by $\Im(m)=\sigma/(2\iu)$.
Moreover, choosing $m$ on this line shows that
\begin{equation}\label{nchi2}
\norm{\varphi(\cdot,\lambda)}_c^2=\frac{\Im(m)}{\Im(\lambda)}=\frac{\sigma(c,\la)}{2\iu\Im(\la)}.
\end{equation}
We denote the half-plane $\{m\in\bb C: \sigma(c,\la)/(2\iu\Im(\la))\leq\Im(m)/\Im(\la)\}$ by $H(c,\la)$.

We shall now allow $c$ to vary.
Suppose $c'<c$, $\norm{\psi(\cdot,\lambda)}_c=0$, and $m\in H(c,\la)$.
Then we have
\[\norm{\varphi(\cdot,\lambda)}_{c'}^2 \leq \norm{\varphi(\cdot,\lambda)}_c^2 \leq \frac{\Im(m)}{\Im(\lambda)}.\]
This means that the point $m$ is also in $H(c',\la)$.
In other words, $H(c,\la)$ is contained in $H(c',\la)$ for $c'\leq c$ and $\Im(\sigma(\cdot,\la))$ is a monotone function.
We define $\sigma(b,\la)=\lim_{c\to b}\sigma(c,\la)$ (which may be $\pm\iu\infty$), when $\norm{\psi(\cdot,\lambda)}_c=0$ for all $c$.

Similarly, if $\norm{\psi(\cdot,\lambda)}_{c'}>0$ and $m\in D(c,\la)$ where $c>c'$, we also have $m\in D(c',\la)$, i.e., $D(c,\la)$ is contained in $D(c',\la)$ for $c'\leq c$.
The radii $r(c,\la)$ are non-increasing as $c$ tends to $b$ and we define $r(b,\la)=\lim_{c\to b}r(c,\la)$.

Finally, if $\norm{\psi(\cdot,\lambda)}_{c'}=0$ but $\norm{\psi(\cdot,\lambda)}_{c}>0$ (implying that $c'<c$) then $D(c,\la)\subset H(c',\la)$.
In this case a half-plane shrinks to a disk as $x$ increases from $c'$ to $c$.
Perhaps it is worthwhile mentioning that on the Riemann sphere the inverse images, under stereographic projection, of both half-planes and disks are disks\footnote{A disk on a sphere is the intersection of a ball with the sphere.} on the sphere.

Now suppose that $\norm{\psi(\cdot,\lambda)}_{c}=0$ for all $c$ implying that the definiteness condition is violated.
Then $H(b,\la)=\bigcap_{0<c<b} H(c,\la)$ is either empty or equal to a half-plane.
On the other hand, if $\norm{\psi(\cdot,\lambda)}_{c}$ is eventually positive, the disks $D(c,\la)$ shrink down to a disk or a singleton, which, in either case, we call $D(b,\la)$.

The first consequence of these considerations is stated in the following lemma.
\begin{lem}\label{L.3.1}
If $\la\not\in\Lambda$, the equation $Ju'+qu=\la wu$ has at least one non-trivial solution in $\mc L^2(w)$.
\end{lem}

\begin{proof}
If $\|\psi(\cdot,\la)\|=0$, then $\psi(\cdot,\la)$ is a non-trivial solution in $\mc L^2(w)$.
Otherwise $\chi_m(\cdot,\la)$ is in $\mc L^2(w)$ when $m$ is chosen in $D(b,\la)$.
\end{proof}

We will now investigate the relationship between the Möbius transform associated with $\la$ and the one associated with $\ol$.
First recall that $\ol$ is not in $\Lambda$ if $\la$ is not in $\Lambda$ since, as mentioned above, $\Lambda$ is symmetric with respect to the real axis.
Let $\mc A=U_{11}(c,\ol)$, $\mc B=U_{21}(c,\ol)$, $\mc C=U_{12}(c,\ol)$, and $\mc D=U_{22}(c,\ol)$.
Then, as $c$ is a point of continuity, Lemma \ref{L_Temp_01} gives
\begin{multline}\label{olla}
\mc C \ov{\mc D}-\ov{\mc C}\mc D=\psi(c,\ol)^*J\psi(c,\ol)=|\tau(c,\ol)|^2\ov{\psi(c,\la)^* J\psi(c,\la)}\\ =|\tau(c,\ol)|^2 (\ov CD -C\ov D).
\end{multline}
Similarly we have
$$\mc A \ov{\mc D}-\mc B\ov{\mc C}
=|\tau(c,\ol)|^2(\ov{A}D-\ov{B}C)
\sas
\mc A \ov{\mc B}-\ov{\mc A}\mc B
=|\tau(c,\ol)|^2(\ov{A}B-A\ov{B}).$$
We also need to consider $\mc A \mc D-\mc B\mc C=\psi(c,\ol)^\top J\varphi(c,\ol)$.
Here we get
$$\mc A \mc D-\mc B\mc C=\tau(c,\ol)^2\ov{\psi(c,\la)^\top J\varphi(c,\la)}=\tau(c,\ol)^2(\ov{AD-BC}).$$
These observations show that the disk $D(c,\ol)$ is the conjugate of the disk $D(c,\la)$ and the half-plane $H(c,\ol)$ is the conjugate of $H(c,\la)$ (whichever the case may be).
In particular, $r(c,\ol)=r(c,\la)$.
For this radius we will now obtain another expression again with the help of Lagrange's identity \eqref{Lagrange}.
Let $(v,g)=(\psi(\cdot,\overline{\lambda}),\overline{\lambda}\psi(\cdot,\overline{\lambda}))$ and $(u,f)=(\phi(\cdot,\lambda),\lambda\phi(\cdot,\lambda))$.
Then
\[\psi(c,\overline{\lambda})^*J\phi(c,\lambda)-\psi(0,\overline{\lambda})^*J\phi(0,\lambda)
 =\int_{[0,c)}(\psi(\cdot,\ol)^*w(\lambda \phi(\cdot,\lambda))-(\ol\psi(\cdot,\overline{\lambda}))^* w \phi(\cdot,\lambda)).\]
Clearly, the integral on the right is zero and $\psi(0,\overline{\lambda})^*J\phi(0,\lambda)=1$.
Hence, using Lemma \ref{L_Temp_01} and equation \eqref{tau},
$$1=\psi(c,\ol)^*J\phi(c,\lambda)=\ov{\tau(c,\ol)} \psi(c,\la)^\top J\phi(c,\lambda)=\frac{AD-BC}{\tau(c,\la)}.$$
This implies
\begin{equation}\label{radii}
r(c,\la)=\frac{|\tau(c,\lambda)|}{2|\Im(\lambda)| \norm{\psi(\cdot,\lambda)}_c^2}
 =\frac{|\tau(c,\ol)|}{2|\Im(\ol)| \norm{\psi(\cdot,\ol)}_c^2} =r(c,\ol).
\end{equation}

We are now in a position to harvest the fruits of our labor.
Our next result resembles Weyl's classical limit-point/limit-circle alternative.

\begin{thm}\label{def}
Suppose that the definiteness condition holds, that $n_+=n_-$ and that $\la\not\in\Lambda$.
Then we have the following dichotomy.
\begin{enumerate}
  \item $D(b,\la)$ has positive radius, if and only if $\norm{\psi(\cdot,\la)}<\infty$.
  In this case all solutions of $Ju'+qu=\la wu$ are in $\mc L^2(w)$, i.e., we have $n_\pm=2$.
  \item $D(b,\la)$ is a singleton, say $\{m_0\}$, if and only if $\norm{\psi(\cdot,\la)}=\infty$.
  In this case we have that $\chi_{m_0}(\cdot,\la)\in\mc L^2(w)$ and $n_\pm=1$.
\end{enumerate}
\end{thm}

\begin{proof}
First note that $\mc L_0=\{0\}$ implies that $\psi(\cdot,\la)$ cannot have $0$ norm and therefore we will eventually have disks rather than half-planes.
Also the two statements are contrapositives of each other, so we need to prove only the forward direction of each.

If $D(b,\la)$ is a proper disk, it contains two points $m_1$ and $m_2$.
Then $\chi_{m_1}(\cdot,\la)$ and $\chi_{m_2}(\cdot,\la)$ are linearly independent elements of $\mc L^2(w)$.
Any other solution of $Ju'+qu=\la wu$, in particular $\psi(\cdot,\la)$, is a linear combination of these two and hence also in $\mc L^2(w)$ and we have $n_\pm=2$.

If $D(b,\la)=\{m_0\}$, then $r(b,\la)=r(b,\ol)=0$.
Assume, by way of contradiction, that $\psi(\cdot,\la)$ has finite norm.
Since $\chi_{m_0}(\cdot,\la)$, which also has finite norm, and $\psi(\cdot,\la)$ are linearly independent, we have $n_\pm=2$ and hence that $\psi(\cdot,\ol)$ also has finite norm.
But now \eqref{radii} implies that $\lim_{c\to b}|\tau(c,\la)|=0$ as well as $\lim_{c\to b}|\tau(c,\ol)|=0$.
This is impossible and proves that $\norm{\psi(\cdot,\la)}=\infty$.
\end{proof}

In the case when the definiteness condition is violated but we still have $n_+=n_-$ the situation is similar to what we just proved.

\begin{thm}\label{indef}
Suppose that the definiteness condition does not hold, but that $n_+=n_-$ and that $\la\not\in\Lambda$.
Then we have the following four cases.
\begin{enumerate}
  \item Assume $\norm{\psi(\cdot,\la)}=0$.
  \begin{enumerate}
    \item $H(b,\la)$ is a half-plane if and only if $\norm{\varphi(\cdot,\la)}<\infty$. In this case we have $n_\pm=1$.
    \item $H(b,\la)$ is empty if and only if $\norm{\varphi(\cdot,\la)}=\infty$. In this case we have $n_\pm=0$.
  \end{enumerate}
  \item Assume $\norm{\psi(\cdot,\la)}>0$.
  \begin{enumerate}
    \item $D(b,\la)$ has positive radius if and only if $\norm{\psi(\cdot,\la)}<\infty$. In this case we have $n_\pm=1$.
    \item $D(b,\la)$ is a singleton, say $\{m_0\}$, if and only if $\norm{\psi(\cdot,\la)}=\infty$. In this case we have $n_\pm=0$.
  \end{enumerate}
\end{enumerate}
\end{thm}

\begin{proof}
First assume that $\norm{\psi(\cdot,\la)}=0$ when we only deal with half-planes.
Our claims in (1) follow from \eqref{nchi2}.

If $\norm{\psi(\cdot,\la)}>0$ we have, eventually, disks to consider.
If $D(b,\la)$ is a proper disk, then all solutions of $Ju'+qu=\la wu$ are in $\mc L^2(w)$, particularly $\psi(\cdot,\la)$.
Since $\dim\mc L_0=1$ we have $n_\pm=1$.

If $D(b,\la)=\{m_0\}$, we repeat the corresponding argument in the proof of Theorem \ref{def} to show that $\norm{\psi(\cdot,\la)}$ is indeed infinite.
The norm of $\chi_{m_0}(\cdot,\la)$ is $0$ and hence $n_\pm=0$.
\end{proof}

We asked in Theorems \ref{def} and \ref{indef} that $n_+=n_-$.
We now investigate under what conditions this assumption fails.

\begin{thm}\label{T.3.4}
$n_+\neq n_-$ if and only if exactly one of $\norm{\psi(\cdot,\la)}$ and $\norm{\psi(\cdot,\ol)}$ is infinite.
\end{thm}

\begin{proof}
First we need to emphasize that $\norm{\psi(\cdot,\la)}=0$ if and only if $\norm{\psi(\cdot,\ol)}=0$.
Indeed, it follows from \eqref{cdcd} and \eqref{olla} that
$$\norm{\psi(\cdot,\ol)}_c^2=\frac{\mc C\ov{\mc D}-\ov{\mc C}\mc D}{2\iu\Im(\ol)}=|\tau(c,\ol)|^2\frac{\ov CD-C\ov D}{2\iu\Im(\ol)} =|\tau(c,\ol)|^2\norm{\psi(\cdot,\la)}_c^2$$
which proves the claim.

Now suppose that $0<\norm{\psi(\cdot,\la)}<\infty$ and $\norm{\psi(\cdot,\ol)}=\infty$.
The latter assumption means that the deficiency index for the half-plane containing $\ol$ is $1-\dim\mc L_0$ and that $D(b,\ol)$ is a singleton.
Hence $D(b,\la)=\{m_0\}$ is also a singleton.
But since $\chi_{m_0}(\cdot,\la)$ and $\psi(\cdot,\la)$ are linearly independent the deficiency index for the half-plane containing $\la$ is $2-\dim\mc L_0$.
This proves one direction of our claim.

To prove the converse, we distinguish three cases.
(1) $\norm{\psi(\cdot,\la)}=\norm{\psi(\cdot,\ol)}=0$, (2) $0<\norm{\psi(\cdot,\la)},\norm{\psi(\cdot,\ol)}<\infty$, and (3) $\norm{\psi(\cdot,\la)}=\norm{\psi(\cdot,\ol)}=\infty$.
In case (1) we have $n_+=n_-=0$ if $H(b,\la)=\emptyset$ and $n_+=n_-=1$ if $H(b,\la)\neq\emptyset$.
In case (2) we obtain from \eqref{radii} that $|\tau(c,\la)|$ and $|\tau(c,\ol)|$ have finite, non-zero limits as $c$ tend to $b$.
Now Lemma \ref{L.2.3} shows that $n_+=n_-$.
In case (3) each of the equations $Ju'+qu=\la wu$ and $Jv'+qv=\ol wv$ has one solution of finite norm and one solution of infinite norm.
This implies that $n_+=1-\dim\mc L_0=n_-$.
\end{proof}

\begin{cor}
If $n_+\neq n_-$, then one of $|\tau(c,\la)|$ and $|\tau(c,\ol)|$ tends to $0$ and the other to infinity as $c$ tends to $b$.
In particular, if $|\tau(c,\la)|$ does not have a limit (in $[0,\infty]$) as $c$ tends to $b$, then $n_+=n_-$.
\end{cor}

\begin{proof}
By Theorem \ref{T.3.4} we may assume, without loss of generality, that $\norm{\psi(\cdot,\la)}<\infty$ and $\norm{\psi(\cdot,\ol)}=\infty$.
Using \eqref{lagforpsi} and \eqref{olla} gives that
$$\frac{\norm{\psi(\cdot,\la)}_c^2}{\norm{\psi(\cdot,\ol)}_c^2}=\frac1{|\tau(c,\ol)|^2}=|\tau(c,\la)|^2$$
tends to $0$ as $c$ tends to $b$.
\end{proof}

We close this section with two examples.

\begin{exm}
This is essentially Example 5.30 of Lesch and Malamud \cite{MR1964480}.
It shows that $n_+$ may indeed be different from $n_-$.
Let $b=\infty$, $\alpha=0$, $q=\sm{0&0\\ 0&0}$ and $w=(1+\frac{a}{x^2+1})\id+\iu J$ with $a\geq0$.
We have
$$U(x,\la)=\e^{\iu\la x}\begin{pmatrix}\cos(\la t(x))&\sin(\la t(x))\\ -\sin(\la t(x))&\cos(\la t(x))\end{pmatrix}$$
where $t(x)=x+a \arctan(x)$.
$\mc L_0$ is spanned by the constant function $(1,-\iu)^\top$ when $a=0$ but for $a>0$ the problem is definite, i.e., $\mc L_0=\{0\}$.
We have $\tau(x,\la)=\exp(2\iu\la x)$ which tends to $0$ when $\Im(\la)>0$ and to infinity when $\Im(\la)<0$.
For $\la$ in the upper half-plane $\psi(\cdot,\la)$ is in $\mc L^2(w)$ and hence all solutions are.
For $\la$ in the lower half-plane only the multiples of $\chi_{-\iu}(\cdot,\la)$ are in $\mc L^2(w)$.
\end{exm}

\begin{exm}
This example shows that we may still have $n_+=n_-$ even though $\tau(c,\la)$ tends to $0$ or $\infty$ as $c$ tends to $b$.
Let $b=\infty$, $\alpha=0$, $q=\sm{0&0\\ 0&0}$ and $w=\sm{4&-\iu\\ \iu&1}$.
Then
$$U(x,\la)=\e^{\iu\la x}\begin{pmatrix}\cos(2\la x)&\frac12\sin(2\la x)\\ -2\sin(2\la x)&\cos(2\la x)\end{pmatrix}$$
and $\norm{\psi(\cdot,\la)}_c^2=(\e^{2c\Im(\la)}-\e^{-6c\Im(\la)})/(8\Im(\la))$.
This tends to infinity as $c$ tends to infinity regardless of whether $\Im(\la)>0$ or $\Im(\la)<0$.
On the other hand $\chi_m(\cdot,\la)$ has finite norm when $m=2\iu\Im(\la)/|\Im(\la)|$.
Hence $n_+=n_-=1$.
\end{exm}

\subsection{Bad points are present for \texorpdfstring{$\la$}{lambda}}\label{bps}
In this section we show that one must avoid to take a $\la\in\Lambda$ when trying to determine the deficiency indices from the behavior of the Weyl disks
since, regardless of the nature of the endpoint $b$, the limit-point situation prevails.

This is most easily made clear by considering some simple examples.
We choose $b=\infty$, $\alpha=0$, $q=\sm{0&2\iu\\ -2\iu&2}\delta_1$ and $w=\sm{2&0\\ 0&0}\delta_1$
where $\delta_1$ is the Dirac measure concentrated at $1$.
Note that $b=\infty$ is a regular endpoint.
Then $B_-(1,2\iu)$ is not invertible but $B_+(1,2\iu)$ is.
Consequently, one may extend any solution on $(0,1)$ to $(0,\infty)$.
However, since $B_-(1,2\iu)$ has only rank $1$ the same is true for $U(c,2\iu)=B_+(1,2\iu)^{-1}B_-(1,2\iu)$ when $c>1$.
Here the first column is $\varphi(c,2\iu)$ and the second is $\psi(c,2\iu)$.
In fact $\varphi(c,2\iu)=-(1+\iu)\psi(c,2\iu)$.
Therefore $m=1+\iu$ gives $\chi_m(c,2\iu)=0$ for $c>1$ and this satisfies the boundary condition $(\cos\beta,\sin\beta)\chi_m(c,2\iu)=0$ for any $\beta$.

Next consider $b=\infty$, $\alpha=0$, $q=\sm{0&-2\iu\\ 2\iu&2}\delta_1$ and $w=\sm{2&0\\ 0&0}\delta_1$.
In this case $B_+(1,2\iu)$ is not invertible but $B_-(1,2\iu)$ is.
Now neither $B_-(1,2\iu)\varphi^-(1,2\iu)$ nor $B_-(1,2\iu)\psi^-(1,2\iu)$ are in the range of $B_+(1,2\iu)$ and it is impossible to define these functions on all of $(0,\infty)$.
But the linear combination $\varphi(\cdot,\la)+m\psi(\cdot,\la)$ where $m=1+\iu$ can be extended from $(0,1)$ to $(0,\infty)$ in infinitely many ways since we may add a solution of $Ju'+qu=2\iu wu$ which is $0$ on $(0,1)$ and satisfies $u^+(1)\in\ker B_+(1,2\iu)$.
One of these will satisfy the boundary condition at $c$ but the value of $m$ is unaffected by this.
Thus again $m$ is simply a point, indeed $m=1+\iu$.

If $B_-(x,\la)$ is always invertible we could also determine the Weyl-Titchmarsh coefficient $m$ in the following way.
Define $\eta(\cdot,\la)$ as the solution of $Ju'+qu=\la wu$ satisfying the initial condition $\eta(c,\la)=(-\sin\beta,\cos\beta)^\top$.
We now ask for which factor $p$ can we satisfy the condition $U(0,\la)\sm{1\\ m}=p \eta(0,\la)$ instead of asking which linear combination $\chi$ of $\varphi$ and $\psi$ satisfies the boundary condition $(\cos\beta,\sin\beta)\chi(c,\la)=0$.
This is equivalent to the requirement
$$\begin{pmatrix}1\\ m  \end{pmatrix}=p\begin{pmatrix}\cos\alpha&\sin\alpha\\ -\sin\alpha&\cos\alpha \end{pmatrix}\eta(0,\la)$$
and yields
\begin{equation}\label{altm}
m=\frac{\eta_2(0,\la)\cos\alpha-\eta_1(0,\la)\sin\alpha}{\eta_2(0,\la)\sin\alpha+\eta_1(0,\la)\cos\alpha}.
\end{equation}
In the above example we get, as expected, $m=1+\iu$.
In fact, if there are no bad points for $\la$, it is easy to see that the Weyl-Titchmarsh coefficient $m$ obtained in Section \ref{gps} is the same as the one given by \eqref{altm}.

We conclude by remarking that it may also happen that $B_-(x,\la)$ and $B_+(x,\la)$ fail to be invertible simultaneously.
In this context the following observation is interesting
\begin{thm}
If the entries of $\Delta_q(x)$ and $\Delta_w(x)$ are real and one of the matrices $B_\pm(x,\la)$ is not invertible, then neither is the other one.
Conversely, if both of the matrices $B_\pm(x,\la)$ fail to be invertible for some non-real $\la$, then the entries of $\Delta_q(x)$ and $\Delta_w(x)$ are real.
\end{thm}

\begin{proof}
This follows immediately from equation \eqref{detbpm} which states
\[\det B_-(x,\la)-\det B_+(x,\la)=2\iu(\Im\Delta_q(x)_{12}-\la \Im\Delta_w(x)_{12}).\qedhere\]
\end{proof}

\section*{Acknowledgment}
We wish to thank Steven Redolfi for fruitful discussions.


\begin{thebibliography}{1}

\bibitem{MR0069338}
Earl~A. Coddington and Norman Levinson.
\newblock {\em Theory of ordinary differential equations}.
\newblock McGraw-Hill Book Company, Inc., New York-Toronto-London, 1955.

\bibitem{MR3046408}
Jonathan Eckhardt, Fritz Gesztesy, Roger Nichols, and Gerald Teschl.
\newblock Weyl-{T}itchmarsh theory for {S}turm-{L}iouville operators with
  distributional potentials.
\newblock {\em Opuscula Math.}, 33(3):467--563, 2013.

\bibitem{MR2145077}
W.~Norrie Everitt.
\newblock Charles {S}turm and the development of {S}turm-{L}iouville theory in
  the years 1900 to 1950.
\newblock In {\em Sturm-{L}iouville theory}, pages 45--74. Birkh\"{a}user,
  Basel, 2005.

\bibitem{MR4047968}
Ahmed Ghatasheh and Rudi Weikard.
\newblock Spectral theory for systems of ordinary differential equations with
  distributional coefficients.
\newblock {\em J. Differential Equations}, 268(6):2752--2801, 2020.

\bibitem{MR1964480}
Matthias Lesch and Mark Malamud.
\newblock On the deficiency indices and self-adjointness of symmetric
  {H}amiltonian systems.
\newblock {\em J. Differential Equations}, 189(2):556--615, 2003.

\bibitem{MR745152}
Jesper L\"utzen.
\newblock Sturm and {L}iouville's work on ordinary linear differential
  equations. {T}he emergence of {S}turm-{L}iouville theory.
\newblock {\em Arch. Hist. Exact Sci.}, 29(4):309--376, 1984.

\bibitem{MR0176151}
E.~C. Titchmarsh.
\newblock {\em Eigenfunction expansions associated with second-order
  differential equations. {P}art {I}}.
\newblock Clarendon Press, Oxford, second edition, 1962.

\bibitem{41.0343.01}
Hermann Weyl.
\newblock {\"Uber gew\"ohnliche Differentialgleichungen mit Singularit\"aten
  und die zu\-geh\"origen Entwicklungen willk\"urlicher Funktionen.}
\newblock {\em Math. Ann.}, 68:220--269, 1910.

\end{thebibliography}
\end{document}